\documentstyle[twoside,12pt]{article} 
\def\Fix{\rm {\hbox{Fix}}}

\newtheorem{theorem}{Theorem} 
\newtheorem{lemma}{Lemma} 
\newtheorem{corollary}{Corollary} 
 
\newtheorem{definition}{Definition} 
 
\newtheorem{remark}{Remark} 
\newtheorem{question}{Question}

\newcommand{\fix}{{\rm Fix\ }}

\date{} 

\title{Nielsen zeta function, 3-manifolds and  asymptotic expansions in Nielsen theory}
\author{Alexander Fel'shtyn\thanks{Part of this work was conducted during  author  stay
 in Max-Planck-Institut f\"ur Mathematik in Bonn}}
\begin{document}   
\bibliography{ref} 
\bibliographystyle{plain}
\maketitle

 \bigskip
 
\setcounter{section}{-1}
 
\section{Introduction}
 Before moving on the results of the paper, we briefly describe the few 
basic notions of Nielsen fixed point theory which will be used.
We assume  $X$ to be a connected, compact
polyhedron and $f:X\rightarrow X$ to be a continuous map.
Let $p:\tilde{X}\rightarrow X$ be the universal cover of $X$
and $\tilde{f}:\tilde{X}\rightarrow \tilde{X}$ a lifting
of $f$, ie. $p\circ\tilde{f}=f\circ p$.
Two liftings $\tilde{f}$ and $\tilde{f}^\prime$ are called
{\sl conjugate} if there is a $\gamma\in\Gamma\cong\pi_1(X)$
such that $\tilde{f}^\prime = \gamma\circ\tilde{f}\circ\gamma^{-1}$.
The subset $p(Fix(\tilde{f}))\subset Fix(f)$ is called
{\sl the fixed point class of $f$ determined by the lifting class $[\tilde{f}]$}.
Two fixed points $x_0$ and $x_1$ of $f$ belong to the same fixed point class iff
 there is a path $c$ from $x_0$ to $x_1$ such that $c \cong f\circ c $ ( homotopy relative
    endpoints). This fact can be considered as an equivalent definition of a non-empty fixed point class.
 Every map $f$  has only finitely many non-empty fixed point classes, each a compact
 subset of $X$. 

A fixed point class is called {\sl essential} if its index is nonzero.
The number of lifting classes of $f$ (and hence the number
of fixed point classes, empty or not) is called the {\sl Reidemeister Number} of $f$,
denoted $R(f)$.
This is a positive integer or infinity.
The number of essential fixed point classes is called the {\sl Nielsen number}
of $f$, denoted by $N(f)$.
 
The Nielsen number is always finite.
$R(f)$ and $N(f)$ are homotopy invariants.
In the category of compact, connected polyhedra the Nielsen number
of a map is, apart from in certain exceptional cases,
 equal to the least number of fixed points
 of maps with the same homotopy type as $f$.

Taking a dynamical point on view,
 we consider the iterates of $f$,
 and we may define several zeta functions connected
 with Nielsen fixed point theory (see \cite{f1,f2,fp}).
The Nielsen zeta function of $f$ is defined
 as power series:
 
\begin{eqnarray*}
 N_f(z)
 & := &
 \exp\left(\sum_{n=1}^\infty \frac{N(f^n)}{n} z^n \right).
\end{eqnarray*}
The Nielsen zeta function  $N_f(z)$ is homotopy invariant.
The function $N_f(z)$ has a positive radius
of convergence which has a sharp estimate in
terms of the topological entropy of the map $f$ \cite{fp}.

We begin the article by proving in Section 1 that 
the Nielsen zeta function is 
a rational function or a radical of a rational function   for orientation preserving 
homeomorphisms of  special Haken
 or  special  Seifert  3-manifolds.

In  Section 2, we obtain 
an asymptotic expansion for 
the number of twisted conjugacy classes or  for the number of 
Nielsen fixed point classes whose norm is 
at most $x$  in the case of pseudo-Anosov homeomorphism of surface.

The author would like to thank  M. Gromov,  Ch. Epstein, R. Hill,  L. Potyagailo, R. Sharp, V.G. Turaev 
for stimulating discussions. In particular, Turaev proposed 
the author  in 1987 the conjecture that the Nielsen zeta function is 
a rational function or a radical of a rational function  for homeomorphisms of  Haken
 or Seifert manifolds.

\section{Nielsen zeta function and homeomorphisms of 3-manifolds}

\subsection{Periodic maps  and homeomorphisms of hyperbolic manifolds}

We prove  in corollary 1 of this subsection  that  the Nielsen zeta function is a radical of rational 
function for any  homeomorphism of a compact hyperbolic 3-manifold. 
Lemma 1  and corollary 1 play also important role in the proof of the
main theorem of this section.

We denote $N(f^n)$ by $N_n$. 
We shall say  that $f:X\rightarrow X$ is a  periodic map of period $m$, 
if $f^m$ is  the identity map $id_X:X\rightarrow X$. Let $ \mu(d), d \in N$, 
be the M\"obius function of number theory. As is known, it is given 
by the following equations:  $\mu(d)=0 $ if $d$ is divisible by a square different from one ; 
$\mu(d)=(-1)^k $ if
 $d$ is not divisible  by a square different from one , where $k$ denotes the number of 
 prime divisors of $d$; $ \mu(1)=1$.

We give  the  proof of the following  key lemma  for the completeness.
 
 \begin{lemma} \cite{fp}
Let $f$ be a  periodic map of least period $m$ of the connected
compact polyhedron $X$ . Then the Nielsen
 zeta function is equal to 
$$
N_f(z)=\prod_{d\mid m}\sqrt[d]{(1-z^d)^{-P(d)}},
$$
 where the product is taken over all divisors $d$ of the period $m$, and $P(d)$ is the integer
$$  P(d) = \sum_{d_1\mid d} \mu(d_1)N_{d\mid d_1} .  $$
\end{lemma}
{\sc Proof}
Since $f^m = id $, for each $j, N_j=N_{m+j}$. If  $(k,m)=1$, then  there exist positive integers $t$
and $q$ such that $kt=mq+1$. So $ (f^k)^t=f^{kt}= f^{mq+1}=f^{mq}f=(f^m)^{q}f =  f$.
Consequently, $  N((f^k)^t)=N(f) $. Let two fixed point $x_0$ and $x_1$ belong 
to the same  fixed point class. Then there exists a path $ \alpha $ from  $x_0$ to $x_1$ 
such that $ \alpha \ast (f\circ\alpha)^{-1} \simeq 0$. 
We have $f(\alpha \ast f\circ\alpha)^{-1})=(f\circ\alpha)\ast (f^2\circ\alpha)^{-1} \simeq 0$ 
and a product 
$ \alpha \ast (f\circ\alpha)^{-1}\ast (f\circ\alpha)\ast (f^2\circ\alpha)^{-1} 
=\alpha \ast (f^2\circ\alpha)^{-1}\simeq 0 $. 
It follows that $ \alpha \ast (f^k\circ\alpha)^{-1} \simeq 0$ is derived by the iteration of this process. 
So $x_0$ and $x_1$ belong to the same  fixed point class of $f^k$. 
If two point belong to the different  fixed point classes $f$, 
then they belong to the different  fixed point classes of $f^k$ . 
So, each essential class( class with nonzero index) for $f$  is an essential class for $f^k$;
in addition , different essential classes for $f$ are different essential classes for $f^k$.
So $ N(f^k)\geq N(f) $. Analogously, $ N(f)=N((f^k)^t) \geq  N(f^k) $. Consequently ,
$ N(f)=N(f^k) $. One can prove completely analogously that $ N_d= N_{di} $, if (i, m/d) =1,
where $d$ is a divisor of $m$. Using these series of equal Nielsen numbers, one can regroup
the terms of the series in the exponential of the Nielsen zeta function so as to get logarithmic
functions by adding and subtracting missing terms with necessary coefficient.
We show how to do this first for period $m=p^l$, where $p$ is a prime number . We have the
following series of equal  Nielsen numbers:
 $$
 N_1=N_k, (k,p^l)=1 (i.e., no\, N_{ip}, N_{ip^2},...
.., N_{ip^l},  i=1,2,3,....),
 $$
 $$ 
 N_p=N_{2p}=N_{3p}=.......=N_{(p-1)p}=N_{(p+1)p}= ... (no \,
N_{ip^2}, N_{ip^3}, ...  , N_{ip^l} ) 
$$
etc.; finally,
$$
N_{p^{l-1}}=N_{2p^{l-1}}=..... (no \, N_{ip^l})
$$
and separately the number $N_{p^l}$.\\
 Further,
\begin{eqnarray*}
\sum_{i=1}^\infty \frac{N_i}{i} z^i & = & \sum_{i=1}^\infty \frac{N_1}{i} z^i 
       +\sum_{i=1}^\infty \frac{(N_p -N_1)}{p}\frac{ {z^p}^i}{i} + \\
                                                    & +  &\sum_{i=1}^\infty \frac{(N_{p^2} -(N_p -N_1)- N_1)}{p^2} \frac{{z^{p^2}}^i}{i} + ...\\ 
 & + & \sum_{i=1}^\infty \frac{(N_{p^l} - ...- (N_p -N_1)- N_1)}{p^l}\frac{ {z^{p^l}}^i}{i}\\
 &=& -N_1\cdot \log (1-z) + \frac{N_1-N_p}{p}\cdot
						    \log (1-z^p) +\\
& + & \frac{N_p-N_{p^2}}{p^2}\cdot \log (1-z^{p^2}) + ...\\
& + & \frac{N_{p^{l-1}}-N_{p^l}}{p^l}\cdot \log (1-z^{p^l}).
\end{eqnarray*}
For an arbitrary period $m$ , we get completely analogously,

\begin{eqnarray*}
N_f(z) & = & \exp\left(\sum_{i=1}^\infty \frac{N(f^i)}{i} z^i \right)\\
           & = & \exp\left(\sum_{d\mid m} \sum _{i=1}^\infty \frac{P(d)}{d}\cdot\frac{{z^d}^i}{i}\right)\\
	   & = & \exp\left(\sum_{d\mid m}\frac{P(d)}{d}\cdot \log (1-z^d)\right)\\
	   & = & \prod_{d\mid m}\sqrt[d]{(1-z^d)^{-P(d)}}
\end{eqnarray*}
where the integers $P(d)$ are calculated recursively by the formula
$$
P(d)= N_d  - \sum_{d_1\mid d; d_1\not=d} P(d_1).
$$
Moreover, if the last formula is rewritten in the form
$$
N_d=\sum_{d_1\mid d}\mu(d_1)\cdot P(d_1)
$$ 
and one uses  the M\"obius Inversion law for real function in number theory, then
$$
P(d)=\sum_{d_1\mid d}\mu(d_1)\cdot N_{d/d_1},
$$
where $\mu(d_1)$ is the M\"obius function in number theory. The lemma is proved.

\begin{corollary}
Let $f: M^n \rightarrow M^n , \, n \geq  3$ be a homeomorphism of a compact hyperbolic manifold. 
Then by Mostow rigidity theorem $f$ is homotopic  to periodic homeomorphism $g$.
So lemma 1  applies and the Nielsen zeta function $N_f(z)$ is equal to
 $$
N_f(z)=N_g(z)= \prod_{d\mid m}\sqrt[d]{(1-z^d)^{-P(d)}},
$$
 where the product is taken over all divisors $d$ of the least period $m$ of $g$,  and $P(d)$ is the integer
$  P(d) = \sum_{d_1\mid d} \mu(d_1)N(g^{d\mid d_1}) .  $
 \end{corollary}

The proof of the following lemma  is based on Thurston's theory of 
homeomorphisms of surfaces \cite{th}.

\begin{lemma}\cite{fp}

The Nielsen zeta function of a  homeomorphism $f$ of a compact  
surface $F$ is either a rational function or the radical of a rational function.
\end{lemma}

{\sc Proof} 
 The case  of a surface with $\chi(F)> 0  $  and case of torus were considered in \cite{fp}.
 If  surface has $\chi(F)=0  $  and $F$ is not
a torus then any homeomorphism is isotopic to periodic one (see \cite{jww}) and Nielsen zeta function
is a radical of rational function by lemma 1 .
In the case of a hyperbolic ($\chi(F)< 0$) surface , according to Thurston's classification
 theorem, the homeomorphism $f$ is isotopic either to a periodic or a pseudo-Anosov, 
or  a reducible homeomorphism. In the first case the assertion of the lemma  follows  from lemma 1.
 If $f$ is an  pseudo-Anosov homeomorphism of a compact surface
 then for each $ n>0, N(f^n)=F(f^n)$ \cite{thur}. Consequently, 
in this case the Nielsen zeta function coincides with the Artin-Mazur zeta function: $N_f(z)=F_f(z)$. 
Since in \cite{fash}  Markov partitions are constructed for a pseudo-Anosov homeomorphism,
 Manning's proof \cite{ma} of the rationality  of the Artin-Mazur zeta function for diffeomorphisms 
satisfying Smale's axiom A carries over to the case of pseudo-Anosov homeomorphisms.
Thus , the Nielsen zeta function $N_f(z)$  is also rational. Now if $f$ is isotopic to a 
reduced homeomorphism $\phi$,
 then there exists a reducing system $S$ of disjoint circles $S_1, S_2, ... , S_m$ on $ int F$
 such that \\
  
  1)  each circle $S_i$ does not bound a disk in $F$;\\
  
  2) $S_i$ is not isotopic to $S_j, i\not=j$;\\
    
  3)  the system of circles $S$ is invariant with respect to $\phi$; \\
  
  4)  the system $S$ has an open $\phi$-invariant tubular neighborhood $\eta(S)$ such that
       each $\phi$ -component $ \Gamma_j $  of the set $ F - \eta(S)$ is mapped into itself by some                                                            
       iterate $\phi^{n_j}, n_j  >0$ of the map $\phi$; here  $\phi^{n_j}$ on $\Gamma_j$ is either a   
       pseudo-Anosov or a
       periodic homeomorphism; \\
  
  5)  each band $ \eta(S_i)$ is mapped into itself by some iterate $\phi^{m_i}, m_i > 0$; here      
       $\phi^{m_i}$
       on $ \eta(S_i)$ is a generalized twist (possibly trivial).\\ 
  
  Since the band $ \eta(S_i)$ is homotopically equivalent to the circle $ S^1$  the Nielsen zeta 
function $N_{\phi^m_i}(z)$ is rational (see \cite{fp}).
  The zeta functions $N_\phi (z)$ and $N_{\phi^m_i}(z)$ are connected on the $\phi$ - component
  $ \Gamma_j $ by the formula $N_\phi (z)=\sqrt[n_j]{ N_{\phi^n_j}(z^{n_j})};$ analogously,
  on the band $ \eta(S_i), N_\phi (z)=\sqrt[m_j]{ N_{\phi^m_j}(z^{m_j})}$.
  The fixed points of $\phi^n$, belonging to different components $ \Gamma_j $
  and bands $ \eta(S_i)$ are nonequivalent \cite{i1}, so the Nielsen number $N(\phi^n)$
  is equal to the sum of the Nielsen numbers $ N(\phi^n / \Gamma_j)$ and $ N(\phi^n / \eta(S_i))$
  of $\phi$ -components and bands . Consequently, by the properties of the exponential,
  the Nielsen zeta function $N_{\phi } (z)=N_f(z)$ is equal to the product of the Nielsen zeta functions 
  of the $\phi$- components $ \Gamma_j $ and the bands $ \eta(S_i)$, i.e. is the radical of a
  rational function.

\subsection{Homeomorphisms of Seifert fibre spaces}

Let $M$ be a compact 3-dimensional Seifert fibre space. That is a space which is foliated by 
simple closed curves, called fibres, such that a fibre $L$ has a neighborhood which is either
a solid Klein bottle or a fibred solid torus $T_r$, where $r$ denotes the number of times 
a fibre near $L$ wraps around $L$. A fibre is regular if it has neighborhoods fibre equivalent
to solid torus $S^1\times D^2$, otherwise it is called a critical fibre. 
There is a natural quotient map $p: M \rightarrow  F$ where $F$ is  a 2-dimensional orbifold,
and hence, topologically a compact surface.
The projection of the critical fibres, which we denote by $S$, consists of a finite set of points
in the interior of $F$ together with a finite subcollection of the boundary components.
 The reader is referred to \cite{sc}
for details, for definitions about 2-dimensional orbifolds, their Euler characteristics
 and other properties of Seifert fibre spaces.
An orbifold is hyperbolic if it has negative Euler characteristic. Hyperbolic orbifold admits hyperbolic
structure with totally geodesic boundary.

Any fibre-preserving homeomorphism $f: M \rightarrow  M$ naturally induces a
relative surface homeomorphism of the pair $(F,S)$, which we will denote by $\hat{f}$.
Recall from \cite{sc} that there is a unique orientable Seifert fibre space with 
orbifold $P(2,2)$ -projective plane with cone singular points of order 2, 2.
We denote this manifold by $M_{P(2, 2)}$.

\begin{lemma}\cite{sc, jww}
Suppose  $M$ is a compact orientable Seifert fiber space which is not
$T^3,  S^1\times D^2, T^2\times I $ or  $ M_{P(2, 2)}$. 
Then there is a Seifert fibration $p: M \rightarrow  F$ , so that
any orientation preserving homeomorphism on $M$ is isotopic to a
fiber preserving homeomorphism with respect to this fibration.
\end{lemma}

Observe that if $M$ and $F$ are both orientable, then $M$ admits a coherent orientation
of all of its fibres and the homeomorphism $f$ either preserves fibre orientation of all
the fibres or it reverses fibre orientation.

\begin{lemma}\cite{ke}
Let $M$ be a compact, orientable aspherical, 3-dimensional Seifert fibre space
such that the quotient orbifold $F$ is orientable and all fibres have neighborhoods
of type $T_1$. Let  $f: M \rightarrow  M$  be a fibre-preserving homeomorphism
inducing  $\hat{f }: F \rightarrow  F$. If  $f$ preserves fibre orientation, then Nielsen number
$N(f)=0$. If $f$ reverses fibre orientation, then $N(f)=2N(\hat{f})$.  
\end{lemma}
{\sc Proof}
First, by small isotopy, arrange that $f$ has a finite number of fibres which are mapped to themselves.
If $f$ preserves fibre orientation a further isotopy, which leaves $\hat{f}$ unchanged, ensures that 
none of these fibres contains a fixed point. Thus, $ \fix(f)=\empty $ and so, $N(f)=0$. 
If $f$ reverses fibre orientation there  are exactly 2 fixed points on each invariant fibre.
Since $M$ is aspherical  lemma 3.2 in \cite{sc} ensures that for any invariant fibre the 2 fixed
points on that fibre are in distinct fixed point classes. On the other hand, the restriction
on the fibre types in the hypothesis implies that a Nielsen path in $F$ can always be lifted to a
Nielsen path path in $M$. In fact, there will be two distinct lifts of each path in $F$.
As a result, $f$ has two fixed point classes covering each fixed point class of $\hat{f}$.
Since the index of a fixed point class of $f$ is the same as that of its projection under $p$ 
the result follows.
\begin{remark}
As one can see lemma   deals with a restricted class of Seifert fibre space.
M. Kelly gave some examples in \cite{ke} which indicate how critical fibres of
Seifert fibre space effect Nielsen classes in the case when $ r > 1$  in $T_r$
and which difficulties arise in this case.

\end{remark}
\begin{definition}
A special Seifert fibre space is a Seifert fiber space such that 
the quotient orbifold $F$ is orientable and all fibres have neighborhoods of type $T_1$. 
\end{definition}

\begin{theorem}
Suppose $M$ is a closed orientable aspherical manifold which is a special Seifert fibre space 
 and $ f: M \rightarrow M $ is an orientation preserving homeomorphism. 
Then the Nielsen zeta function $N_f(z)$ is a rational function or a radical of  a rational function.
\end{theorem}

{\sc Proof}
If  $f$ preserves fibre orientation, then $f^n$ also preserves fibre orientation and so, by lemma 4
Nielsen numbers $N(f^n)=0$ for all $n$ and  the Nielsen zeta function  $N_f(z)=1$.
If $f$ reverses fibre orientation , then $f^2$ preserves fibre orientation, $f^3$ reverses fibre orientation
and so on.  Thus, we have  $N(f^{2k+1})=2N({\hat f}^{2k+1})$ and  $N(f^{2k})=0$ for
$k=0, 1, 2, ...$   .  As result, the Nielsen zeta function  $N_f(z)$  equals
\begin{eqnarray*}
N_f(z) & = & \exp\left(\sum_{n=1}^\infty \frac{N(f^n)}{n} z^n \right)\\
&=& \exp\left(\sum_{k=0}^\infty \frac{N(f^{2k+1})}{2k+1} z^{2k+1} \right)\\
&=& \exp\left(\sum_{k=0}^\infty \frac{2N({\hat f}^{2k+1})}{2k+1} z^{2k+1} \right)\\
&=& \exp \left(2\cdot\sum_{n=1}^\infty \frac{N(\hat f^{n})}{n} z^n
 - 2\cdot\sum_{k=1}^\infty \frac{N(\hat f^{2k})}{2k} z^{2k }\right)\\
&=&  ({N_{\hat f}(z)})^2 \cdot \exp -\left(\sum_{k=1}^\infty \frac{N(({{\hat f}^2})^{k})}{k} (z^2)^k \right)\\
&=& (N_{\hat f}(z))^2 /N_{{\hat f}^2}(z^2 ).
\end{eqnarray*}

From  lemma  2   it follows that  Nielsen zeta functions $N_{\hat f}(z)$ and $N_{{\hat f}^2}(z^2 )$ are either
rational functions or the radicals of rational functions.
Consequently,  the Nielsen zeta function $N_f(z)$ is a rational function or a radical of  a rational function.

\subsection{ Main theorem}
Basic concepts  about 3-manifold topology can be found in \cite{jaco}, in particular the
Jaco-Shalen-Johannson  decomposition of Haken manifolds and Seifert fiber spaces.
See \cite{th} for discussions on hyperbolic 3-manifolds.
We recall some basic facts about compact connected orientable 3-manifolds.
A 3-manifold $M$ is irreducible if every embedded 2-sphere bounds an embedded 3-disk.
By the sphere theorem \cite{hem}, an irreducible 3-manifolds is a $K(\pi, 1)$ Eilenberg-MacLane space 
if and only if it is 
a 3-disk or has infinite fundamental group.

A properly-embedded orientable connected surface in a 3-manifold is incompressible
if it is not a 2-sphere and the inclusion induces a injection on the fundamental groups.
An irreducible 3-manifold is Haken if it contains an embedded orientable incompressible surface.
 We use notation $N(X)$ to denote a regular neighborhood of  set $X$.
A  Haken  3-manifold  $M$ can be decomposed
along a canonical set $\cal T$ of  incompressible tori into pieces  such that each component
of $M - N(\cal T)$ is either hyperbolic manifold, or twisted $I$-bundle over Klein bottle, or a Seifert
fiber space with hyperbolic orbifold. 
This decomposition is called  Jaco-Shalen-Johannson decomposition. 
\begin{definition}
 A special Haken manifold is a Haken manifold  $M$ such that each component
of $M - N(\cal T)$ in JSJ decomposition  of $M$ is either hyperbolic manifold, or a twisted $I$-bundle
over Klein bottle,  or a  aspherical special  Seifert fiber space with hyperbolic orbifold.  
\end{definition} 
We need some definitions from the paper \cite{jww}.
\begin{definition} Suppose $f : M \rightarrow M$ is a map, and $A,B$ are $f$-invariant
sets of $M$. If there is a path $\gamma$ from $A$ to $B$ such that
$\gamma \sim f\circ \gamma$ rel $(A,B)$, then we say that $A, B$ are $f$-related.
\end{definition}
\begin{definition}
Suppose $M$ is a compact 3-manifold with torus boundaries.
A map $f : M \rightarrow M$  is standard on boundary if for any component $T$ of
$\partial M$, the map $f/T$ is one of the following types: (1) a fixed point free map;
(2) a periodic map with isolated fixed points; (3) a fiber preserving, fiber orientation reversing map with respect to some 
$S^1$ fibration of $T$.
\end{definition}
\begin{definition}
A map $f$ on a compact 3-manifold $M$ is said to have $FR$-property(fixed-point relating property)
if the following is true: if $A\in \Fix(f)$ and $B$ is either a fixed point of $f$ or an $f$-invariant
component of $\partial M$, and $A, B$ are $f$-related by a path $\gamma$, then $\gamma$ is
$A,B$ homotopic to a path in $\Fix(f)$.
\end{definition}
\begin{definition}
A map $f$ is a type I standard map if (1) $f$ has $FR$-property, (2) $f$ is standard on boundary,
(3) $\Fix(f)$ consists of isolated points, and (4) $f$ is of flipped pseudo-anosov type at each fixed point.
A map $f$ is type II standard map if it satisfies (1), (2) above, as well as (3) $\Fix(f)$ is a properly
embedded 1-dimensional submanifold and (4) $f$ preserves a normal structure on $\Fix(f)$.
\end{definition}
The main result of this section is the following theorem.
\begin{theorem}
Suppose $M$ is a closed orientable  manifold which is  special Haken manifold, 
and $ f: M \rightarrow M $ is an orientation preserving homeomorphism. 
Then the Nielsen zeta function $N_f(z)$ is a rational function or a radical from rational function.
\end{theorem}

{\sc Proof}
Let $\cal T$ be (possibly empty) set
of invariant tori of the JSJ decomposition of $M$ . Then each component
of $M - Int N(\cal T)$ is either hyperbolic manifold,  or a twisted $I$-bundle over Klein bottle, or
a aspherical special Seifert fiber space with hyperbolic orbifold. Isotop $f$ so that it maps
$N(\cal T)$ homeomorphically to itself.
Suppose that $P$ is a Seifert fibered component of $M - Int N(\cal T)$  such that
$f(P)=P$. By \cite{sc}, $f/P$ is isotopic to a fiber preserving map.
Recall that a torus $T$ in $M$ is a vertical torus if it is union of fibers in $M$.
By \cite{jww}  lemma 1.10, we can find a set of vertical tori $\cal T^*$ in $P$,
cutting $P$ into pieces which are either a twisted $I$-bundle over Klein
bottle, or have  hyperbolic orbifold, and a fiber preserving isotopy of $f$,
so that after isotopy, the restriction of $f$ on each invariant piece has
periodic or pseudo-Anosov orbifold map. 
Adding all such $\cal T^*$ to $\cal T$, we get a collection of tori
$\cal T'$, such that :  (1) $ f(N(\cal T')= N(\cal T')$; (2) each component
$M_i$ of $M - Int N(\cal T')$  either ist a twisted $I$-bundle over Klein
bottle or has hyperbolic orbifold; (3) If $f$ maps a Seifert fibered component
$M_i$ to itself, and if $M_i$ has hyperbolic orbifold, then $f$ is fiber 
preserving, and the orbifold map $\hat f$ on $F_i$ is either periodic
or pseudo-Anosov.
By \cite{jww} we can isotop $f$ so that  restriction  $f /M_i $ is a standard map
for all $M_i$  and after that we can further isotop $f$, rel $\partial N_j$, on each
component $N_j$ of $N(\cal T')$, so that it is a standard map on $N_j$.
By the definition of standard maps, $\Fix(f)$ intersects each of $M_i$ and $N_j$
in points and  $1$-manifolds, so $\Fix(f)$  is a disjoint union of points, arcs, and
circles. 
  Jiang, Wang and Wu proved \cite{jww}theorem 9.1 that  different components of $\Fix(f)$  are not
equivalent in Nielsen sense. We repeat their proof here  for the completeness.
If exist two different components $C_0, C_1$ of $\Fix(f)$   which are equivalent, then there is 
a path $\alpha $ connecting  $C_0, C_1$  such that $f\circ\alpha \sim \alpha $ rel $\partial$.
Denote  by $T''$ the set of tori $\partial N(\cal T')$.  Among all such $\alpha$, choose
one such that  the number of components $\sharp (\alpha - \cal T'')$  of the set 
$\alpha - \cal T''$ is minimal. In below we will find another such curve $\alpha ''$ with
$\sharp (\alpha'' - \cal T'') <\sharp (\alpha - \cal T'')$, which would  contradict
the choice of $\alpha$. 
   Let $D$ be a disk, and let $h: D\rightarrow M$ be a homotopy 
$f\circ \alpha \sim \alpha $ rel $\partial$. We may assume  that $h$
is transverse to $\cal T''$, and $\sharp h^{-1}(\cal T'')$ is minimal
among all such $h$. Then $h^{-1}(\cal T'')$  consists of a properly 
embedded $1$-manifold on $D$, together with possibly one or two 
isolated points mapped to the ends of $\alpha$. $\cal T''$ is $\pi_1$-injective
in $M$, so one can modify $h$ to remove all circles in $h^{-1}(\cal T'')$.
Note that $h^{-1}(\cal T'')$  must contain some arcs, otherwise $\alpha$
would lie in some $M_i$ or $N(T_j)$, which is impossible because the restriction of $f$ 
in each piece has FR-property.
Now consider an outermost arc $b$ in $h^{-1}(\cal T'')$  . Let $\beta=h(b)$.
The ends of $\beta$ can not both be on $\alpha$, otherwise 
we can use the outermost disk to homotope $\alpha$ and reduce $\sharp  ( \alpha  - \cal T'')$,
contradicting the choice of $\alpha$. Since $f$ is a homeomorphism, the same
 thing is true for $f\circ\alpha$. Therefore, $\beta$ has one end on each of $\alpha$ and
$f\cdot\alpha$.  The arc $b$  cutt off a disk $\Delta$ on $D$ whose interior is disjoint
from $h^{-1}(\cal T'')$. The boundary of $\Delta$ gives rise to a loop
$h(\partial \Delta)=\alpha_1\cup\beta\cup(f\circ\alpha_1)^{-1}$, where
$\alpha_1$ is subpath of $\alpha$ starting from an end point $x$ of $\alpha$.
let $T$ be the torus in $\cal T''$ which contains $\beta$. Then the restriction
of $h$ on $\Delta$ gives a homotopy $\alpha_1\sim f\cdot\alpha_1$ 
rel$(x,T)$. Since $f$ has FR-property on each component of $M - Int N(\cal T)$ 
and $N(\cal T)$, by definition there is a path $\gamma$ in $\Fix (f)$ such that 
$\gamma \sim \alpha_1$ rel $(x,T)$. Since $\gamma$ is in $\Fix(f)$, the path
$\alpha '=\gamma ^{-1}\cdot \alpha$ has the property that
$f\circ\alpha '=(f\circ\gamma ^{-1})\cdot (f\circ\alpha)=
\gamma ^{-1}\cdot (f\circ\alpha) \sim \gamma ^{-1}\cdot\alpha =\alpha '$
rel $\partial$. Since $\alpha_1 \sim \gamma$ rel $(x,T)$, the path 
$\gamma ^{-1}\cdot \alpha_1$ is rel $\partial$ homotopic to a path $\delta$ on $T$.
Write $\alpha=\alpha_1\cdot\alpha_2$. Then $\alpha '= \gamma ^{-1}\cdot \alpha=
(\gamma ^{-1}\cdot \alpha_1)\cdot\alpha_2 \sim \delta\cdot\alpha_2$ rel $\partial $.
By a small perturbation on $\delta$, we get a path $\alpha ''\sim \alpha '$ rel
$\partial $ such that $\sharp (\alpha'' - \cal T'') <\sharp (\alpha - \cal T'')$.
Since $f\circ\alpha ''\sim f\circ\alpha '\sim \alpha '\sim \alpha ''$ rel $\partial$,
this contradicts the minimality of  $\sharp (\alpha - \cal T'')$.

We can prove in the same way that the fixed points of  an each iteration $f^n$, 
belonging to different components $M_i$ and $N_j$ are nonequivalent , so the Nielsen number $N(f^n)$
  equals  to the sum of the Nielsen numbers $ N(f^n / M_i)$ and $ N(f^n / N_j)$
  of components  . Consequently, by the properties of the exponential,
  the Nielsen zeta function $N_{f } (z)$ is equal to the product of the Nielsen zeta functions for
the induced homeomorphisms of the components $ M_i$ and  $N_j$. 
By  corollary 1   the Nielsen zeta function of the hyperbolic
component $M_i$ is a radical of a rational function.  The Nielsen zeta function 
is a radical of a rational function for the component $M_i$ which is a aspherical special
Seifert fibre space by the theorem  1  .
The Nielsen zeta function is a rational function or a radical of a rational function 
for a homeomorphism of torus  or Klein bottle by lemma 2  .
This  implies that the Nielsen zeta function is a rational function or a radical of a rational function 
for  the induced homeomorphism  of component $N_j$ which is $I$-bundle over torus
 and for the induced homeomorphism of component $M_i$ which is twisted I-bundle over Klein bottle.
So, the Nielsen zeta function  $N_{f } (z)$ of the homeomorphism $f$  of the 
whole manifold $M$ is a rational function 
or a radical of a rational function 
as a  product of the Nielsen zeta functions for
the induced homeomorphisms of the components $ M_i$ and  $N_j$.

\section{ Asymptotic  expansions for fixed point classes and twisted conjugacy classes 
of pseudo-Anosov homeomorphism }

\subsection{Twisted conjugacy and Reidemeister Numbers.}

Let $\Gamma$ be a group and $\phi:\Gamma\rightarrow \Gamma$ an endomorphism.
Two elements $\alpha,\alpha^\prime\in \Gamma$ are said to be
$\phi-conjugate$ iff there exists $\gamma \in \Gamma$ with
$$
\alpha^\prime=\gamma  \alpha   \phi(\gamma)^{-1}.
$$
We shall write $\{x\}_\phi$ for the $\phi$-conjugacy class
 of the element $x\in \Gamma$.
The number of $\phi$-conjugacy classes is called the $Reidemeister$
$number$ of $\phi$, denoted by $R(\phi)$. If $\phi$ is the identity map then 
the $\phi$-conjugacy classes are the usual conjugacy classes in the group $\Gamma$. 

In \cite{for} we have conjectured that the Reidemeister number to be infinite
as long as the endomorphism is injective and  group has exponential growth.
Below we prove this conjecture for surface groups and pseudo-Anosov maps
and,  in fact,  we obtain  an asymptotic expansion for 
the number of twisted conjugacy classes  whose norm is at most $x$.

Let $f:X\rightarrow X$ be given, and let a
specific lifting $\tilde{f}:\tilde{X}\rightarrow\tilde{X}$ be chosen
as reference.
Let $\Gamma$ be the group of
 covering translations of $\tilde{X}$ over $X$.
Then every lifting of $f$ can be written uniquely
as $\gamma\circ \tilde{f}$, with $\gamma\in\Gamma$.
So elements of $\Gamma$ serve as coordinates of
liftings with respect to the reference $\tilde{f}$.
Now for every $\gamma\in\Gamma$ the composition $\tilde{f}\circ\gamma$
is a lifting of $f$ so there is a unique $\gamma^\prime\in\Gamma$
such that $\gamma^\prime\circ\tilde{f}=\tilde{f}\circ\gamma$.
This correspondence $\gamma\rightarrow\gamma^\prime$ is determined by
the reference $\tilde{f}$, and is obviously a homomorphism.
 
\begin{definition}
The endomorphism $\tilde{f}_*:\Gamma\rightarrow\Gamma$ determined
by the lifting $\tilde{f}$ of $f$ is defined by
$$
  \tilde{f}_*(\gamma)\circ\tilde{f} = \tilde{f}\circ\gamma.
$$
\end{definition}
 
It is well known that $\Gamma\cong\pi_1(X)$.
We shall identify $\pi=\pi_1(X,x_0)$ and $\Gamma$ in the following way.
Pick base points $x_0\in X$ and $\tilde{x}_0\in p^{-1}(x_0)\subset \tilde{X}$
once and for all.
 Now points of $\tilde{X}$ are in 1-1 correspondence with homotopy classes of paths
in $X$ which start at $x_0$:
for $\tilde{x}\in\tilde{X}$ take any path in $\tilde{X}$ from $\tilde{x}_0$ to 
$\tilde{x}$ and project it onto $X$;
conversely for a path $c$ starting at $x_0$, lift it to a path in $\tilde{X}$
which starts at $\tilde{x}_0$, and then take its endpoint.
In this way, we identify a point of $\tilde{X}$ with
a path class $<c>$ in $X$ starting from $x_0$. Under this identification,
$\tilde{x}_0=<e>$ is the unit element in $\pi_1(X,x_0)$.
The action of the loop class $\alpha = <a>\in\pi_1(X,x_0)$ on $\tilde{X}$
is then given by
$$
\alpha = <a> : <c>\rightarrow \alpha \cdot  c = <a\cdot c>.
$$
Now we have the following relationship between $\tilde{f}_*:\pi\rightarrow\pi$
and
$$
f_*  :  \pi_1(X,x_0) \longrightarrow \pi_1(X,f(x_0)).
$$

\begin{lemma}
Suppose $\tilde{f}(\tilde{x}_0) = <w>$.
Then the following diagram commutes:
$$
\begin{array}{ccc}
  \pi_1(X,x_0)  &  \stackrel{f_*}{\longrightarrow}  &  \pi_1(X,f(x_0))  \\
                &  \tilde{f}_* \searrow \;\; &  \downarrow w_*   \\
                &                           &  \pi_1(X,x_0)
\end{array}
$$
where $w_*$ is the isomorphism induced by the path $w$.
\end{lemma}
In other words, for every  $\alpha = <a>\in\pi_1(X,x_0)$ , we have
$$
\tilde{f}_*(<a>)=<w(f \circ a)w^{-1}>
$$ 
\begin{remark}
In particular, if $ x_0 \in p(Fix( \tilde f))$ and $ \tilde x_0 \in Fix( \tilde f)$,then $ \tilde{f}_*=f_*$.
\end{remark} 
 
\begin{lemma}
Lifting classes of $f$( and hence  fixed point classes, empty or not)
 are in 1-1 correspondence with
$\tilde{f}_*$-conjugacy classes in $\pi$,
the lifting class $[\gamma\circ\tilde{f}]$ corresponding
to the $\tilde{f}_*$-conjugacy class of $\gamma$.
We therefore have $R(f) = R(\tilde{f}_*)$.
\end{lemma}
 
We shall say that the fixed point class $p(Fix(\gamma\circ\tilde{f}))$,
which is labeled with the lifting class $[\gamma\circ\tilde{f}]$,
{\it corresponds} to the $\tilde{f}_*$-conjugacy class of $\gamma$.
Thus $\tilde{f}_*$-conjugacy classes in $\pi$ serve as coordinates for fixed
point classes of $f$, once a reference lifting $\tilde{f}$ is chosen.

\subsection{Asymptotic expansions}

We assume $X$ to be a compact surface of negative Euler characteristic and
 $f:X\rightarrow X$ is a pseudo-Anosov homeomorphism, i.e. there is a number $\lambda >1$ 
and a pair of transverse measured foliations $(F^s,\mu^s)$ and $(F^u,\mu^u)$ 
such that $f(F^s,\mu^s)=(F^s,\frac{1}{\lambda}\mu^s)$ and $f(F^u,\mu^u)=(F^u,\lambda\mu^u)$.
The mapping torus $T_f$ of $f:X\rightarrow X$ is the space obtained from $X\times [0 ,1]$ 
 by identifiing $(x,1)$ with $(f(x),0)$ for all $x\in X$. It is often more convenient to regard $T_f$ 
as the space obtained from $X\times [o,\infty )$ by identifiing $(x,s+1)$ with $(f(x),s)$ for 
all $x\in X,s\in [0 ,\infty )$. On $T_f$ there is a natural semi-flow 
$\phi :T_f\times [0,\infty )\rightarrow T_f, \phi_t(x,s)=(x,s+t)$ for all $t\geq 0$.
Then the map  $f:X\rightarrow X$ is the return map of the semi-flow $\phi $. 
A point $x\in X$ and a positive number $\tau >0$ determine the orbit 
curve $\phi _{(x,\tau )}:={\phi_t(x)}_{0\leq t \leq \tau}$ in $T_f$. 
The fixed points and periodic points of $f$ then correspond to closed orbits of various periods.
Take the base point $x_0$ of $X$ as the base point of $T_f$. According to van Kampen Theorem 
the fundamental group $G :=\pi_ 1(T_f,x_0)$ is obtained from $\pi $ by adding a new 
generator $z$ and adding the relations $z^{-1}gz=\tilde f_*(g)$ for all $g\in \pi =\pi _1(X,x_0)$, 
where $z$ is the generator of $\pi_1(S^1,x_o)$
This means that $G$ is a semi-direct product $G=\pi \bowtie  Z$ of $\pi$ with $Z$.

We now describe some known results  

\begin{lemma}
If $\Gamma$ is a group and $\phi$
 is an endomorphism of $\Gamma$ then
 an element $x\in \Gamma$ is always
 $\phi$-conjugate to its image $\phi(x)$.
\end{lemma}

{\it Proof.}
If $\gamma=x^{-1}$ then one has immediately
 $\gamma x = \phi(x) \phi(\gamma)$.
The existence of a $\gamma$ satisfying this equation implies
 that $x$ and $\phi(x)$ are $\phi$-conjugate.

\begin{lemma}
Two elements $x,y$ of $\pi$ are $\tilde f_*$-conjugate iff $xz$ and $yz$ are conjugate 
in the usual sense in $G$. Therefore $R(f) = R(\tilde{f}_*)$ is the number of usual conjugacy 
classes in the coset $\pi \cdot z$ of $\pi$ in $G$. 
 \end{lemma}

{\it Proof.}
if $x$ and $y$ are $\tilde{f}_*$-conjugate then there is a $\gamma \in \pi$ such 
that $\gamma x=y\tilde{f}_*(\gamma)$. This implies $\gamma x=yz\gamma z^{-1}$ 
and therefore $\gamma(xz)=(yz)\gamma$ so $xz$ and $yz$ are conjugate in the usual sense in $G$. 
Conversely suppose $xz$ and $yz$ are conjugate in $G$. 
Then there is a $\gamma z^n \in G$ with $\gamma z^n xz=yz\gamma z^n$. 
From the relation $zxz^{-1}=\tilde{f}_*(x)$ we 
obtain $\gamma \tilde{f}_*^n(x)z^{n+1}=y\tilde{f}_*(\gamma) z^{n+1} $ 
 and therefore $\gamma \tilde{f}_*^n(x)=y\tilde{f}_*(\gamma)$. 
This shows that $\tilde{f}_*^n(x)$ and $y$ are $\tilde{f}_*$-conjugate. 
However by lemma 7  $x$ and $\tilde{f}_*^n(x)$ are $\tilde{f}_*$-conjugate, so $x$ and $y$ 
must be $\tilde{f}_*$-conjugate.

There is a canonical projection
 $\tau : T_f \to R/Z$ given by
 $(x,s)\mapsto s$.
This induces a map
 $\pi_1(\tau):G=\pi_1(T_f,x_0)\to Z$.

We see that the Reidemeister number $R(f)$
 is equal to the number of homotopy
 classes of closed paths $\gamma$ in $T_f$
 whose projections onto $R/Z$
 are homotopic to the path
 $$
\begin{array}{cccl}
 \sigma :&[0,1] & \to     & $ $ \!\!\! R/Z \\
         &  s   & \mapsto & s.
\end{array}
 $$

Corresponding to this there is a group theoretical interpretation of
 $R(f)$ as the number of
 usual conjugacy classes of elements $\gamma\in\pi_1(T_f)$
 satisfying $\pi_1(\tau)(\gamma) = z$.

 \begin{lemma}\cite{thu, otal}
Interior of  the mapping torus $ Int(T_f)$ admits a hyperbolic structure of finite volume if and only if
$f$ is isotopic to pseudo-Anosov homeomorphism.
\end{lemma}

So, if the surface $X$ is closed and $f$ is isotopic to pseudo-Anosov homeomorphism
the mapping torus $T_f$ can be realised as a hyperbolic 3-manifold ,  $H^3/G$, where $H^3$
is the Poincare upper half space $\{(x,y,z): z  > 0, (x,y) \in R^2\}$ with the metric
$ds^2=(dx^2 +dy^2 + dz^2)/z^2$. 
The closed geodesics on a hyperbolic manifold are in one-to-one correspondence with
the free homotopy classes of loops. These classes of loops are in one-to-one correspondence with the 
conjugacy classes of loxodromic elements in the fundamental group of the hyperbolic manifold .
This correspondences allow Ch. Epstein ( see \cite{eps}, p.127) to study the asymptotics of such functions 
as   $p_n(x)$= \# \{ primitive closed geodesics of lenght less than $x$ represented by an element of the
form $gz^n$ \}  using the Selberg trace formula. A primitive closed geodesic is one which is not 
an iterate of another closed geodesic. Later, Phillips and Sarnak \cite{ps} generalised results of Epstein 
 and obtained  for $n$-dimensional hyperbolic manifold the asymptotic of the number 
of primitive closed geodesics of length at most $x$ lying in fixed homology class. 
The proof of this result makes  routine use of the Selberg trace formula.  In the more general case
of variable negative curvature such asymptotic was obtained by Pollicott and Sharp \cite{pos}.
They used dynamical approach based on the geodesic flow.
We will only need an asymptotic   for $p_1(x)$.  Note, that closed geodesics represented by an
element of the form $gz$ are  automatically primitive, becouse they wrap 
exactly once around the mapping torus (once around generator $z$). 
We have following asymptotic expansion \cite{eps, ps, pos} 
\begin{equation}
p_1(x)=  \frac{e^{hx}}{x^{3/2}}( \sum_{n=0}^{N}\frac{C_n}{x^{n/2}} + O(\frac{1}{x^{N/2}}) ),
\end{equation}
for any $N>0$,
where $h=\dim T_f  - 1=2$ is the topological entropy of the geodesic flow on the unit-tangent bundle 
$ST_f$, the constant  $C_0 > 0$ depend on the volume of  hyperbolic 3-manifold $T_f$.
N. Anantharaman \cite{an} has shown that constants $C_n$ vanish if $n$ is odd.
So, we have following asymptotics
\begin{equation}
p_1(x) \sim  C_0 \frac{e^{hx}}{x^{3/2}},   x \rightarrow \infty
\end{equation}

Notation.  We write $f(x)=g(x) +O(h(x)) $ if there exists $C>0$ such that $\mid f(x)-g(x) \mid \leq Ch(x) $.
We write $f(x)\sim g(x)$  if  $\frac{f(x)}{g(x)}\rightarrow 1 $  as  $ x \rightarrow \infty$

Now, using  one-to-one correspondences  in lemma 6  and   lemma  8  
we  define  a norm of fixed point class, or corresponding to him 
lifting class, or corresponding to them twisted 
conjugacy class  $\{g\}_{\tilde f_*}$ in the fundamental group
of surface $\pi =\pi _1(X, x_0)$  as the length of the primitive closed geodesic  $\gamma $  on  $T_f$,
 which represented by an element of the form $gz$. So, for example, the norm function $l^*$ 
on the set of twisted  conjugacy classes equals  $l^*=l\circ B$,
where $l$ is  length function on geodesics($l(\gamma)$ is the length of the 
primitive closed geodesic  $\gamma$ ) and $B$ is bijection
between the set of twisted 
conjugacy class  $\{g\}_{\tilde f_*}$ in the fundamental group
of surface $\pi =\pi _1(X, x_0)$  and the set of closed geodesics represented by an
elements of the form $gz$ in the fundamental group $G :=\pi_ 1(T_f,x_0)$ .
We introduce following  counting functions\\ 
 FPC(x)= \# \{fixed point classes of $f$ of norm  less than $x$ \},\\ 
L(x)= \# \{ lifting  classes of $f$ of norm  less than $x$ \},\\
Tw(x)= \#\{ twisted conjugacy  classes for  $\tilde f_*$  in the fundamental group
of surface  of norm  less than $x$ \} 

\begin{theorem}
Let $X$ be a closed surface of negative Euler characteristic and $f:X\rightarrow X$ 
is a pseudo-Anosov homeomorphism.Then 
$$ FPC(x)=L(x)=Tw(x)= 
 \frac{e^{2x}}{x^{3/2}}( \sum_{n=0}^{N}\frac{C_n}{x^{n/2}}  + O(\frac{1}{x^{N/2}}) ),$$
where the constant  $C_0 > 0$ 
depend on the volume of  hyperbolic 3-manifold $T_f$,
constants $C_n$ vanish if $n$ is odd.
\end{theorem}

{\sc Proof}
The proof follows from lemmas 8 and 9  and   asymptotic expansion (1).
\begin{corollary}
For pseudo-Anosov homeomorphism of surface  the Reidemeister
number is infinite.
\end{corollary}
We can generalise Theorem 3 in following way
\begin{theorem} Let $M$ be a compact manifold  and 
$f:M\rightarrow M$  is a homeomorphism. Suppose that the mapping torus
$T_f$ admits a Riemannian metric of negative sectional curvatures  . Then there exist constants
$C_0, C_1, C_2, ...$ with $C_0>0$ such that
$$ FPC(x)=L(x)=Tw(x)= 
 \frac{e^{hx}}{x^{3/2}}( \sum_{n=0}^{N}\frac{C_n}{x^{n/2}} + O(\frac{1}{x^{N/2}}) ) $$
for any $N>0$,
where $h>o$ is the topological entropy of the geodesic flow on the unit-tangent bundle 
$ST_f$  and constants $C_n$ vanish if $n$ is odd.
\end{theorem}

\begin{question}
For which  compact manifolds $M$ and homeomorphisms $f:M\rightarrow M$ 
the mapping torus
$T_f$ admits a Riemannian metric of negative curvature or negative sectional curvatures ? 
\end{question}

\begin{question}
How to define the norm of the fixed point class or twisted conjugacy class
in general case? Is it exist twisted Selberg trace formula for the discrete group?
Such formula can give asymptotic expansion of counting function 
for twisted conjugacy classes.

\end{question}

\bigskip
 
 Institut f\"ur Mathematik, 

Ernst-Moritz-Arndt-Universit\"at Greifswald 
 
Jahn-strasse 15a, D-17489 Greifswald, Germany.
 
{\it E-mail address}: felshtyn@mail.uni-greifswald.de
\bigskip

\end{document}